\documentclass[12pt, reqno]{amsart}
\usepackage{amscd,amsthm,amsfonts,amssymb,amsmath,euscript}
\usepackage[dvips]{graphicx}
\usepackage[dvips]{graphics}
\usepackage[matrix,arrow]{xy}

\theoremstyle{definition}

\theoremstyle{remark}

\newcommand{\ZZ}{{\mathbb Z}}

\newcommand{\PP}{{\mathbb P}}
\newcommand{\CC}{{\mathbb C}}

\newcommand{\FF}{{\mathbb F}}
\newcommand{\Aff}{{\mathbb A}}

\newcommand{\Lem}{{\bf Lemma }}

\newcommand{\Rem}{{\bf Remark}}

\newcommand{\Prop}{{\bf Proposition}}
\newcommand{\Proof}{{\bf Proof}}
\newcommand{\Comment}{\bf Commentary}

\newfont{\smallskob}{cmbx7 scaled\magstep4}
\newfont{\bigskob}{cmbx12 scaled\magstep4}

\newcommand{\pic}{\mathrm{Pic}\,}

\newcommand{\tit}{Two non-conjugate embeddings of $S_3\times \mathbb Z_2$ into
the Cremona group II}

\begin{document}

\begin{title}
\tit
\end{title}

\begin{abstract}
We prove more precisely the following main result of~\cite{Isk1}.
There is two non-conjugate embeddings of $S_3\times \ZZ_2$ into
the Cremona group that is given by the linear action of this group
on the plane and the action on the two-dimensional torus.
\end{abstract}

\author{V.\,A.\,Iskovskikh}

\thanks{The work was partially supported by RFFI grant No. 050100353.}

%\address{Steklov Institute of Mathematics, 8 Gubkin street, Moscow 117966, Russia} %

\address{Steklov Institute of Mathematics, 8 Gubkin street, Moscow 117966, Russia} %

\email{iskovsk@mi.ras.ru}

\maketitle

\section{Introduction and the main result}

\subsection{}
\label{question}
Remind the main result which is inspired
by~\cite{LPR}. Consider two following actions of the group
$G\simeq S_3\times \ZZ_2$ on the rational surfaces $P$ and $T$
(where $S_3$ is the symmetric group and $\ZZ_2$ is the cyclic
group of order two).

\begin{description}
  \item[(I) $P$ is the plane $x+y+z=0$] $S_3$ acts by the permutations
  on the coordinates and $\ZZ_2$ is given by $(x,y,z) \rightarrow
  (-x,-y,-z)$;
  \item[(II) $T$ is the two-dimensional torus $xyz=1$] $S_3$ acts by the permutations on the coordinates and
  $\ZZ_2$ is given by $(x,y,z) \rightarrow
  (x^{1-},y^{-1},z^{-1})$.
\end{description}

In the both cases we have embedding of this group into
two-dimensional Cremona group $Cr_2(\CC)$. The question is: are
the images of these embeddings of $G$ conjugate in $Cr_2(\CC)$ or
not?

Our answer is the following.

\subsubsection \Prop \label{Proposition:main_result} {\it The
images of these embeddings of $G$ are not conjugate in
$Cr_2(\CC)$. In the other words, there is not exist
$G$-equivariant birational map between $G$-surfaces $P$ and $T$. }

\subsection \Comment Together with~\cite{LPR} this result means
that the simple algebraic group of type $\mathbb G_2$ is not
Cayley group in the following sense. Let $G$ as in~\cite{LPR} be
the connected algebraic group (we let it only in this commentary).
Let $\mathfrak g$ be its Lie algebra. Consider $G$ as the
algebraic variety $X$ with the action of $G$ on itself by the
conjugation $gxg^{-1}$, $g\in G$, $x\in X$, and the Lie algebra
$\mathfrak g$ as the algebraic variety $Y$ with adjoint action
$Ad_Gg(y)=gyg^{-1}$, $g\in G$, $y\in Y$. The group $G$ is called
\emph{Cayley group} if there exist $G$-equivariant birational map
$\lambda\colon X \dashrightarrow Y$, i. e.
$\lambda(gxg^{-1})=Ad_Gg(\lambda(x))$. The question about Cayley
property is reduced to some property of character lattice
$\hat{T}$ of the maximal torus $T\in G$ with the natural action of
Weyl group $W$ on $\hat{T}$. Many canonical linear algebraic
groups are proved to be non-Cayley ones. However, there also exist
a lot of exceptions.

The maximal torus of a group of type $\mathbb G_2$ is
two-dimensional, and its Weil group is $W\simeq S_3\times \ZZ_2$.
The question about Cayley property is exactly
problem~\ref{question}. Our result proves that $\mathbb G_2$ is
not Cayley. It is surprising that in~\cite{LPR} proved that the
group $\mathbb G_2\times G_m^2$ (where $G_m$ is the multiplicative
group of the field with the trivial action of Weyl group) is a
Cayley group. Thus, $\mathbb G_2$ is stable Cayley in this sense.
It would be interesting to know, is the group $\mathbb G_2\times
G_m$ Cayley or not?

In our work~\cite{Isk1} we give the sketchy proof of
proposition~\ref{Proposition:main_result} based on the general
method of factorization of $G$-equivariant birational maps between
rational $G$-surfaces (see~\cite{Ma1},~\cite{Isk3}). $G$ here is a
finite group. However we give complete and more precise proof of
it because of great interest in this fact.

\subsection{}
\label{1.3} As in~\cite{Isk1} we use the conception of rational
$G$-surface for proof. Here $G\subset Cr_2(\CC)$ is the finite
subgroup in the Cremona group (see, for
instance,~\cite{Ma1},~\cite{Isk2}). In our case $G\simeq S_3\times
\ZZ_2$ and it acts on the compactifications of $Y$ and $X$ of
surfaces $P$ and $T$ (which correspond to cases (I) and (II)
respectively). Describe these actions more precisely.

\medskip

{\bf The case (I).} The surface $Y=\PP^2$. Let $(u_0,u_1,u_2)$ be
the homogeneous coordinates on $\PP^2$ and $x=\frac{u_1}{u_0}$,
$y=\frac{u_2}{u_0}$ and $z=-\frac{u_1+u_2}{u_0}$. Then the action
of $G$ on $Y$ is given in the matrix form as follows (we use the
homogenous coordinates).

$$
    \mbox{The involution }\sigma_{xy}=(xy)=
    \left(\begin{array}{ccc}
         1 &     0 &     0 \\
         0 &     0 &     1 \\
         0 &     1 &     0 \\
    \end{array}\right);
$$

$$
    \mbox{the 3-cycle }\sigma_{xyz}=(xyz)=
    \left(\begin{array}{rrr}
         1 &     0 &     0 \\
         0 &     0 &     1 \\
         0 &    -1 &    -1 \\
    \end{array}\right);
$$

$$
    \mbox{the involution }\tau=
    \left(\begin{array}{rrr}
         1 &     0 &     0 \\
         0 &    -1 &     0 \\
         0 &     0 &    -1 \\
    \end{array}\right) \mbox{is the generator of $\ZZ_2$.}
$$
Fixed elements: the point $(1,0,0)$.

Invariant elements: the line $L_0=(u_0=0)$. The group $\ZZ_2$ acts
on $L_0$ trivially and $S_3$ acts as the standard irreducible
two-dimensional linear representation. A unique $G$-invariant
$0$-orbit of length $2$ is $\{(0,1,\lambda),(0,\lambda,1)\}$,
where $\lambda=e^{\frac{2\pi i}{3}}$ is the cubic root of $1$.
There are also two $G$-invariant $0$-orbits of length $3$, which
are given by $\{(0,0,1),(0,1,-1),(0,1,0) \}$ and
$\{(0,1,1),(0,1,-2),(0,-1,1) \}$.

There is $G$-equivariant pencil of lines with $G$-fixed point
$(1,0,0)$ and $G$-invariant section (the line $L_0$). The group
$S_3$ acts on the base of this pencil and on the section $L_0$.
The group $\ZZ_2$ acts on the fibers by $t  \rightarrow -t$.

\medskip

{\bf The case (II).} The surface $X$ is the most natural smooth
compactification of two-dimensional torus $T\colon (xyz=1)$.
Consider the threefold $\PP^1\times \PP^1\times\PP^1$ with
homogeneous coordinates $(x_1,x_0)\times(y_1,y_0)\times(z_1,z_0)$
and coordinates $(x=\frac{x_1}{x_0}, y=\frac{y_1}{y_0},
z=\frac{z_1}{z_0})$. Then this compactification is given by
\begin{equation}
\label{1}
x_1y_1z_1=x_0y_0z_0.
\end{equation}
Under the Segre map $\PP^1\times \PP^1\times \PP^1 \hookrightarrow
\PP^7$ the equation $(1)$ is the $G$-invariant equation of the
hyperplane $\PP^6\supset X$. The surface $X$ is given there by
the quadratic equations that are restrictions on $\PP^6$ of the
equations which give the Segre map. Under this embedding $X$ is a
smooth projective del Pezzo of degree $6$ in $\PP^6$.

$G$ acts in (\ref{1}) as follows. $S_3$ acts by the coordinate
permutations induced by the permutations of factors of
$\PP^1\times\PP^1\times\PP^1$ (for example $\sigma_{xy}\colon
(x_1,x_0) \longleftrightarrow (y_1,y_0)$). $\ZZ_2$ acts as
follows.
$$
(x_1,y_1,z_1) \longleftrightarrow (x_0,y_0,z_0).
$$
These actions are linear under the Segre map.

According to the classification of $G$-minimal rational
$G$-surfaces (see~\cite{Isk2},~\cite{Mo}), $X\subset \PP^6$ is the
minimal del Pezzo $G$-surface of degree $6$ with $G$-invariant
Picard group $\pic^G(X)=\ZZ(-K_X)$, where $K_X$ is the canonical
class. The configuration of $(-1)$-curves on $X$ is a unique
$G$-orbit and have the form of regular $6$-angle. Its equation is
the infinite hyperplane section
\begin{equation}\label{2}
  x_0y_0z_0=0.
\end{equation}
These lines are one-dimensional fibers of three birational
projections $X\rightarrow \PP^1\times \PP^1$, two one-dimensional
fibers in each of them.

\medskip

Further we need other $G$-birational equivalent to $X$ models
$X_0$, $X_1$, $X_2$ and classification of some $0$-dimensional (of
small length) and $1$-dimensional orbits on them.

\subsection{}
\label{1.4}
Start from the surface $X$ given by (\ref{1}). $\pic
X$ is generated by $(-1)$-curves and all $(-1)$-curves forms a
unique $G$-orbit, which is equivalent to $-K_X$, so
$\pic^G_X=\ZZ(-K_X)$. We may consider an affine model $T\colon
(xyz=1)$ for studying of $G$-orbits, because the orbit on infinite
is known (it is all $(-1)$-curves) and there is not $0$-orbit of
length less than $6$.

\subsubsection \Lem
\label{1.4.1} {\it Let $\mathcal A\subset T$ be the
$G$-equivariant $0$-orbit of length $d<6$. Then $\mathcal A$ is one of the
following.
\begin{description}
  \item[$d=1$] a unique fixed point $P=(1,1,1)$;
  \item[$d=2$] a unique orbit $\{ P_1=(\lambda,\lambda,\lambda),P_{-1}=(\lambda^{-1},\lambda^{-1},\lambda^{-1})\}$,
  where $\lambda=e^{\frac{2\pi i}{3}}$ is a cubic root of $1$ {\rm (}notice that $\lambda^{-1}=\lambda^2${\rm )};
  \item[$d=3$] a unique orbit $\{ Q_1=(1, -1, -1), Q_2=(-1,1,-1), Q_3=(-1,-1,1)$;
  \item[$d=4$] there is no such orbits;
  \item[$d=5$] there is no such orbits, because $5\nmid 12=|G|$.
\end{description}
}

Indeed, one can see this from the equation $xyz=1$ and the action
of $G$ on $T$. In the case $d=4$ the stabilizer of the point is
$\langle\sigma_{xyz}\rangle=\ZZ_3$ acting by $(x,y,z)\rightarrow
(y,z,x)$, which means that $x=y=z$.

\subsubsection{}
Now find some of $G$-orbits, which consist of rational curves.

\medskip

$\mathcal T=\Gamma_x+\Gamma_y+\Gamma_z\sim -K_X$. This $G$-orbit
consists of the curves of genus $0$ which contain the fixed point
$P=(1,1,1)$ and is given by $x=1$, $y=1$, $z=1$. In the projective
form these equations on $\PP^1\times\PP^1$ are as follows.
%\begin{multline}
$$
\Gamma_x\colon y_1z_1-y_0z_0=0,\\
$$
\begin{equation}\label{3}
\Gamma_y\colon x_1z_1-x_0z_0=0,\\
\end{equation}
$$
\Gamma_z\colon x_1y_1-x_0y_0=0.\\
$$
%\end{multline}

Each of these curves $\Gamma$ is smooth one of genus $0$ with
$-K_X\Gamma =2$, $\Gamma^2=0$, and each two of them intersect by a
unique point $P$. Blow up this point. Then after $G$-equivariant
contraction of these curves we have three $G$-conjugate points on
the image of exceptional $(-1)$-curve of blow-up on the model
$X_2\subset \PP^3$, see below (this is example of $G$-equivariant
link $\Phi_{6,1}$, see section $2$).

This birational $G$-map $\gamma\colon X\dashrightarrow X_2$ is a
projection from the tangent plane to the fixed point $P$ with the
embedding $X\subset \PP^6$ described above. Under this projection
the point $P$ blows up, and three conics passing through $P$
contract.

\medskip

$\mathcal D=\Delta_x+\Delta_y+\Delta_z\sim -2K_X$. This is a triple
of smooth curves of genus $0$ with equations $y=z$, $z=x$, and
$x=y$ respectively. In the projective form we have
$$
\Delta_x\colon x_1y_1^2-x_0y_0^2=0,\\
$$
\begin{equation}\label{4}
\Delta_y\colon y_1z_1^2-y_0z_0^2=0,\\
\end{equation}
$$
\Delta_z\colon z_1x_1^2-z_0x_0^2=0.\\
$$
All these curves $\Delta$ intersect in three points $P$, $P_1$ and
$P_{-1}$ with $-K_X\Delta=4$ and $\Delta^2=2$.

\subsubsection{}
\label{subsec143}
 Curves $\Gamma$ and $\Delta$ are components of
the fibers of $G$-invariant pencil of rational curves
$\Pi=|-K_X-2P-P_1-P_{-1}|$ (i. e. the pencil of curves from the
linear system $|-K_X|$ which contain points $P_1$, $P_{-1}$ and
twice point $P$). Moreover, $\Gamma_x+\Delta_x\sim
\Gamma_y+\Delta_y\sim \Gamma_z+\Delta_z$ are the fibers of this
pencil. Equations (\ref{3}) and (\ref{4}) involve that
$\Gamma_x\cap \Delta_x=\{P,Q_1\}$, $\Gamma_y\cap
\Delta_y=\{P,Q_2\}$, $\Gamma_z\cap \Delta_z=\{P,Q_3\}$.

\medskip

$\mathcal E=E_x+E_y+E_z\sim -K_X$. This triple of smooth curves of
genus $0$ is given by $x=-1$, $y=-1$, $z=-1$ accordingly. One can
see from the equations that $E_x\ni Q_2, Q_3$, $E_y\ni Q_1,Q_3$,
$E_z\ni Q_1,Q_2$, and $Q_i$ are the only intersection points for
the components of $\mathcal E$. As above, $-K_XE=2$, $E^2=0$, so
after blow up of $Q_1$, $Q_2$ and $Q_3$ strict transforms
$E_x^\prime$, $E_y^\prime$, $E_z^\prime$ are $G$-conjugate triple
of $(-2)$-curves which do not intersect each other. They may be
$G$-equivariant contracted to three ordinary double points on the
cubic surface in $\PP^3$.

\subsection{} Consider now $G$-birational model $X_0$ of
$G$-surface $X$, namely, standard projectivisation of torus
$T\subset \Aff^3\subset \PP^3$
\begin{equation}\label{5}
X_0\colon xyz=w^3,
\end{equation}
where $(x,y,z,w)$ are homogenous coordinates on $\PP^3$. This
model differs from $X$ only on infinity $w=0$. There are not $6$
lines as on $X$ but three ones, which are given by $x=w=0$,
$y=w=0$ and $z=w=0$. There are $3$ ordinary double points
$(0,0,1,0)$, $(0,1,0,0)$, and $(1,0,0,0)$. $G$ acts on $X_0$ not
linear (as on $X\subset \PP^6$ does). More precisely, $S_3$ acts
linearly by permutations of $x,y,z$ and $\ZZ_2$ acts by birational
involution $(x,y,z,w)\rightarrowtail
(x^{-1},y^{-1},z^{-1},w^{-1})$.

However there is another cubic model $X_1\subset \PP^3$ such that
$G$ acts linearly on it. There are three ordinary double points as
on $X_0$. This model is given by $G$-equivariant birational
projection
$$
p_1\colon X\subset \PP^6 \dashrightarrow \PP^3\supset X_1
$$
from the plane $\langle Q_1, Q_2, Q_3\rangle$, which is a linear
span of $Q_1,Q_2,Q_3\in X$ in $\PP^6$. As we notice in~\ref{1.4},
conics $E_x, E_y,E_z$ contract to the singular points.

Choose a homogenous coordinates $(x,y,z,w)$ in $\PP^3$ such that
the image of $P$ is a point $(0,0,0,1)$ and $\ZZ_2$ acts by
involution $(x,y,z,w)\rightarrowtail (-x,-y,-z,w)$. Then $S_3$
acts by coordinate permutations as on $X_0$. Under the projection
$p_1\colon X\dashrightarrow X_1$ three singular points (the images
of $E_x$, $E_y$, and $E_z$) lie on the infinity $w=0$. The images
of conics $\Gamma_x$, $\Gamma_y$, and $\Gamma_z$ are three lines
passing through $P_0=p_1(P)$. Thus, $P_0$ is the Eckardt point,
whose tangent plane is $G$-invariant, i. e. it is given by linear
equation $x+y+z=0$.

So, the equation of $X_1\subset \PP^3$ is the following.
\begin{equation}\label{6}
X_1\colon xyz-aw^2(x+y+z)=0, \ \ \ a\in \CC.
\end{equation}

We can put $a=\frac{1}{3}$ for convenience. Then the singular
points are $(1,0,0,0)$, $(0,1,0,0)$, $(0,0,1,0)$, and the images
of $P_1$ and $P_{-1}$ are $(1,1,1,1)$ and $(-1,-1,-1,-1)$. Three
lines $x=0$, $y=0$, $z=0$ lie in the tangent plane $x+y+z=0$.

\subsection \Rem As at each cubic surface, there is Geiser's birational
involution concerned with $P_0$, i. e. the projection from this
point to the plane composed with the automorphism of the double
covering. Obviously this birational involution is $G$-equivariant
and in our case it is biregular because $P_0$ is the Eckardt
point. $G$-equivariant birational Bertini involution is concerned
with two $G$-conjugate points $P_1, P_{-1}\in X$. Its action on
$\pic X_1$ differs from the standard one
because $X_1$ is singular and the points $P_1, P_{-1}$ are
not in the general position. This means that the third
intersection point of the line $\langle P_1, P_{-1}\rangle \subset
\PP^3$ and $X_1$ (i.e. $P_0$) lies on the $(-1)$-curves.

\subsection{} There is also one $G$-equivariant birational model
of $X$. Namely, quadric $X_2\subset \PP^3$ of type
\begin{equation}\label{7}
X_2\colon xy+yz+zx=3w^2.
\end{equation}

It is given from $X_1$ by $G$-equivariant birational transform
$(x,y,z,w)\rightarrowtail (x^{-1},y^{-1},z^{-1},w^{-1})$. Indeed,
this coordinate change with the multiplication to the common
factor transform the equation (\ref{6}) with $a=\frac{1}{3}$ to
the equation (\ref{7}). The action of $G$ on $(x,y,z,w)$ are still
linear. Moreover, $X_2$ is nothing but the image of $X\subset
\PP^6$ under the linear projection
$$
p_2\colon X\subset \PP^6 \dashrightarrow \PP^3\supset X_2
$$
from the tangent plane to $X$ at the $G$-fixed point $P\in X$. The
projection $p_2$ blows up $P$ to the conic $C_0\colon (w=0)\subset
X_2$ and contracts three conics $\Gamma_x, \Gamma_y,
\Gamma_z\subset X$ to the $G$-invariant $0$-orbit $\mathcal
A=\{(1,0,0,0), (0,1,0,0), (0,0,1,0)\}\subset C_0$.

The curve $C_0$ is the image of the point $P_0\in X_1$ and the
$0$-orbit $\mathcal A$ is the image of the lines $x=0, y=0, z=0$
under our map $X_1 \dashrightarrow X_2$.

\subsection \Lem {\it
\begin{description}
  \item[(i)] The variety $X_2$ is smooth and there is no $G$-fixed points on it;
  \item[(ii)]Conics $C_0=(w=0)$ and $C_1=(x+y+z=0)$ are $G$-invariant and on $C_0$ acts only $S_3$ by a unique
two-dimensional irreducible representation and $G$ acts effective
on $C_1$;
  \item[(iii)] $G$-invariant $0$-orbits of length $d<6$ on
  $X_2$ are:
\begin{description}
  \item[$d=2$] $\{P_1=(1,1,1,1), P_{-1}=(-1,-1,-1.1)\}$, the images on $X_2$ of $P_1$ and $P_{-1}$ on X {\rm (}or $X_1${\rm )};
$\{R_1=(1,\lambda, \lambda^{-1}),
R_2=(\lambda,1,\lambda^{-1})\}=C_0\cap C_1$, the basic points of
the conic pencil on $X_2$ generated by $C_0$ and $C_1$,
$\lambda=e^{\frac{2\pi i}{3}}$
  \item[$d=3$] $\mathcal A=\{(1,0,0,0),(0,1,0,0),(0,0,1,0)\},\ \ \ \ \ \ \ \ \ \ \ \ \ \ \ \ \ \ \ \ \ \ \ \ \ \ \ \ \ \ \ $
               $\mathcal B=\{(-1,2,2,0),(2,-1,2,0),(2,2,-1,0)\}$,
               $\mathcal A, \mathcal B\subset C_0$;
  \item[$d=4$] there is no such orbits;
  \item[$d=5$] there is no such orbits;
\end{description}
    \item[(iv)] let
$$
\Pi_0=(t_0(x+y+z)+t_1w=0), \ \ (t_0,t_1)\in \PP^1,
$$
be the pencil of conics generated by $C_0$ and $C_1$; then on the
base of this pencil $\PP^1$ acts only $\ZZ_2=G/S_3\colon
(t_0,t_1)\rightarrowtail (-t_0,t_1)$ with two fixed points $(1,0)$
and $(0,1)$; there is two reducible fibers $F_1$ and $F_{-1}$ on
$\Pi_0$ with intersection points $P_1$ and $P_{-1}$ of the
components with coordinates $(1,-3)$ and $(1,3)$ accordingly; on
the irreducible fibers $G$ acts without fixed points;
    \item[v] there is another pencil of conics
$$
\Pi_1=\{u_0(x-y)+u_1(y-z)+u_2(z-x)=0; \ \ u_0+u_1+u_2=0, \ \ (u_0,
u_1, u_2)\in \PP^2\}
$$
passing through $P_1, P_{-1}\in X_2$; it is an image of $\Pi$ from clause
\ref{subsec143} under the projection $p_2\colon X
\dashrightarrow X_2$; on the base of this pencil
$\PP^1=(u_0+u_1+u_2=0)$ as usual acts $S_3$ by an two-dimensional
irreducible representation, $\ZZ_2$ acts on each fiber with fixed
points $\Pi_1\cap C_0$ and invariant $2$-section $\Pi_1\cap C_1$.

There is two reducible fibers $G_1$ and $G_2$ on $\Pi_1$ with
intersection points of the components $R_1$ and $R_2$ with
coordinates $(\lambda^{-1},1,\lambda)$ and
$(\lambda^{-1},\lambda,1)$ on the base accordingly.
\end{description}
}

\medskip

All these statements can be straightforward check using the
equations and the action of $G$.

\section{Proof of the main result}

\subsection{} The idea of the proof of
proposition~\ref{Proposition:main_result} is the same as in
the~\cite{Isk1}. We should show that there is no $G$-equivariant
birational map between $G$-surfaces  $X$ and $Y=\PP^2$, defined in
section \ref{1.3}. According to the general theory
(see~\cite{Isk3}), each $G$-equivariant birational map can be
decomposed (by algorithm given by $G$-equivariant Sarkisov
program) to the sequence of $G$-equivariant elementary birational
maps--links. The classification of such links is obtained in the
algebraic case if $G$ is finite Galois group, acting on smooth
$k$-minimal rational surfaces over the perfect field $k$
(see~\cite{Isk3}).

However, according to the general conception of rational
$G$-surfaces (see~\cite{Ma1},~\cite{Isk2}), this classification
may be applied to the geometrical case, i. e. to the case of
action of the finite group $G$ on smooth projective minimal
rational $G$-surfaces. Their classification is also known
(see~\cite{Isk2} and $G$-equivariant $2$-dimensional Mori
theory~\cite{Mo}).

Some differences between geometrical and algebraical cases
concerned with the description of $G$-orbits. For example, in the
geometrical case (with non-trivial $G$-action) $G$-orbits are
closed proper subset. In the algebraical case, for instance, for
del Pezzo surfaces of great degrees, the existence of one
$k$-point (orbit of degree $1$) involves everywhere density of
such $k$-points in the Zariski topology.

The theorem of decomposition of any birational $G$-map to a
sequence of links is known (see~\cite{Isk3}, $2.5$) and all links
are determined by their centers, i. e. $G$-orbits of dimension
$0$, all points of which are in the general position on the
corresponding $G$-surface. So, for proof of the main result we
need the following.

\begin{itemize}
\item[{\bf a)}] Choose from the general classification
(see~\cite{Isk3}, $2.6$) the links which is concerned with our
case. Indicate minimal $G$-surfaces on which acts these links (as
$G$-equivariant birational maps), starting from $X$.
\item[{\bf b)}] Classify all zero-dimensional $G$-orbits of length
$d< \deg X_i=K_{X_i}^2$ on all of such surfaces $X_i$ and choose
those of them, whose points are in the general position (this is
necessary condition for the existence of links).
\end{itemize}

Remind, that the points $x_1,x_2, \ldots, x_d\in X$ are in the
general position, if

\begin{itemize}
  \item[{\bf 1)}] $X$ is del Pezzo surface, $\sigma\colon X^\prime\rightarrow X$ is the simultaneous blow up of all
  $x_1, \ldots, x_d$ and $X^\prime$ is also del Pezzo
  surface, or
  \item[{\bf 2)}] $\pi\colon X\rightarrow \PP^1$ is a conic bundle, then none of this points lye on the
  degenerated fibers and at most one point of $x_1,\ldots, x_d$ lies on the fiber; this is necessary and sufficient
  condition of the existence of link with center in $x_1,\ldots,x_d$, i. e. the birational transform
$$
\xymatrix{ %X&\ar@{-->}[rr] & X^\prime
%&&W\ar@{->}[ld]_{g}\ar@{->}[rd]^{f}&&\\%
&X\ar@{->}[rd]\ar@{-->}[rr]&&X^\prime,\ar@{->}[ld]& \\
&&\PP^1
%&&\PP^1\ar@{->}[ld]_{g}\ar@{->}[rd]^{f}&&
}
$$
which is blow up of $x_1,\ldots,x_d$ and a contraction of
strict transforms of the fibers on which they lie.
\end{itemize}

\subsection{} Suppose that there exists $G$-equivariant birational
map $\chi\colon X \dashrightarrow Y$. Then $\chi$ can be
decomposed to the composition of links. According to the
particular algorithm of decomposition and considering all
possibilities for links we will see that there is no link whose
image is $\PP^2$. So we will obtain a contradiction with the
decomposition theorem. So there is no such $\chi$, which means
that two embeddings $G$ to the Cremona group that we define in (I)
and (II), section 1.1, are not conjugate in $Cr_2(\CC)$.

By construction from the algorithm choose the very ample
$G$-invariant linear system $\mathcal H^\prime=\mathcal H_Y$ on
$Y$, for example, $\mathcal H^\prime=|-K_Y|$. Let $\mathcal
H=\mathcal H_X=\chi^{-1}_*\mathcal H^\prime$ be its strict
transform on $X$. By the definition of strict transform the linear
system $\mathcal H$ has no fixed components and $\dim \mathcal
H=\dim \mathcal H^\prime$ ($\mathcal H$ has basic points if $\chi$
is not morphism). $\pic^G(X)=\ZZ(-K_X)$, so $\mathcal H\sim
-aK_X$, $a\in \ZZ_{>0}$.

The map $\chi\colon X\rightarrow Y$ is not isomorphism, so, by
$G$-equivariant Noether inequality (see~\cite{Isk3}, $2.4$),
linear system $\mathcal H$ has maximal singularity, which is a
zero-dimensional $G$-orbit $x\subset X$ of length $d=d(x)$ and
multiplicity $mult_x\mathcal H=r=r(x)>a$ at all points $x_i\in x$.
We have $\mathcal H^2=a^2 K_X^2>r^2d$ ($\mathcal H$ is mobile
linear system without fixed components, passes through all points
of $G$-orbits $x$ with multiplicity $r$, and determines a
birational map) so $d<K_X^2$, in our case $d<K_X^2=6$. The length
of the orbit is divided by the order of the group. So, the maximal
singularities by which we can construct links can be only
$G$-orbits of length $d=1$, $2$, $3$ or $4$. We find such orbits
in lemma~\ref{1.4.1}.

\subsection{} Now, for every such $G$-orbit we should check
are its points in the general position or not and if yes, then
choose from the classification (\cite{Isk3}, $2.6$) appropriate
link $\Phi\colon X\dashrightarrow X_1$ (do not confuse with $X_1$
from paragraph $1.5$). Algorithm of decomposition is the
following. The map $\chi_1\colon X_1\dashrightarrow Y$ in the
composition $\chi=\chi_1\circ \Phi$ is given by the linear system
$\mathcal H_1=\Phi_*(\mathcal H)\sim -a_1K_{X_1}$ with $a_1<a$. We
have $a\in \ZZ_{>0}$, so the decomposition should be terminated on
the isomorphism after the finite number of steps. Start to find
links. There is no $G$-orbit of length $4$ on $X$
(see~\ref{1.4.1}), so start from $d=3$.

{\bf The case $d=3$.} There is one such orbit $(Q_1,Q_2,Q_3)$.
Show that in this case there is no such ``untwisting'' link
$\Phi$. Indeed, the points $Q_1$, $Q_2$, $Q_3$ on $X$ are not in
the general position. On the blow up of these three points
$\sigma\colon X^\prime \rightarrow X$ strict transforms
$E^\prime_x$, $E^\prime_y$, $E^\prime_z$ of conics $E_x$, $E_y$,
$E_z$ (see $1.4.2$) are $G$-invariant triple of $(-2)$-curves, so
$-K_X$ is not ample ($-K_X\cdot E^\prime_x=0$).

(As we see in paragraph $1.5$, the projection $p_1\colon
X\dashrightarrow X_1$ from the plane $\langle Q_1, Q_2,
Q_3\rangle\subset \PP^6 \dashrightarrow \PP^3$ blows up these
three points and contracts $E_x$, $E_y$, $E_z$ to the ordinary
double points on $X_1$.)

Notice that in~\cite{Isk1} we mistakenly say that there is such
link $\Phi_{6,3}$. The mistake is the following. We do not check
the condition of generality (though this is not affect to the
final result).

\subsubsection \Rem The matter why the link $\Phi_{6,3}$ does not
exists is easy.  For existence of link with center in the maximal
singularity $x\subset X$ it is necessary and sufficient that under
the blow up $\sigma\colon X^\prime \rightarrow X$of the cycle $x$
on the surface $X^\prime$ there is two extremal rays, one of which
is the exceptional divisor. $\pic^G (X^\prime)=\ZZ\oplus\ZZ$, so the
Mori cone is generated by these extremal rays and $-K_{X^\prime}$
is ample by Kleiman's criterion. The link $\Phi$ is nothing but
the blow up $\sigma$ and the extremal contraction of this second
ray.

In the classical terms the non-generality of position of $Q_1$,
$Q_2$ and $Q_3$ means that this $G$-orbit could not be a maximal
singularity in the linear system $\mathcal H$, because if $r>a$,
then $E_x$, $E_y$ and $E_z$ are fixed components of $\mathcal
H$. Indeed, $\mathcal H\subset |-aK_X-rQ_1-rQ_2-rQ_3|$, the cycles
$E_x-Q_2-Q_3$, $E_y-Q_1-Q_3$, $E_z-Q_1-Q_2$ are effective and must
have non-negative intersection with $\mathcal H$. But $\mathcal
H\cdot (E_x-Q_2-Q_3)=(-aK_x-rQ_1-rQ_2-rQ_3)\cdot
(E_x-Q_2-Q_3)=2a-2r<0$. So we have a contradiction.

\subsubsection{\bf The case $d=2$.}
There is a unique such orbit $x=\{P_1, P_{-1}\}$
(see~\ref{1.4.1}). The pair $P_1$, $P_{-1}$ is in the general
position on $X$. Indeed, each of these points does not lie on the
$(-1)$-curves (which lie on the infinity $x_0y_0z_0=0$), and there
is no curve of genus $0$ and degree $2$ (with respect to $-K_X$)
on which both of these points lie.

There is corresponding link $\Phi=\Phi_{6,2}\colon
X\dashrightarrow X_1$ (see~\cite{Isk3}, $2.6$). Find out what is
$X_1$. Let $\sigma\colon X^\prime\rightarrow X$ be the blow up of
$x=\{P_1,P_{-1}\}$. Then $-K_{X^\prime}$ is ample and
$K_{X^\prime}^2=4$. The linear system $|-K_{X^\prime}|$ gives
$G$-equivariant embedding $X^\prime \subset \PP^4$ into the del
Pezzo surface of degree $4$ (i. e. into an intersection of two
quadrics). The image of the exceptional $(-1)$-curves are two skew
lines whose linear span is a hyperplane $\PP^3$. This image is
$G$-invariant and intersects $X^\prime\subset \PP^4$ by another
$G$-invariant pair of skew lines which is reducible curve of genus
$1$ together with the pair that we blow up.

Contract this residual pair of lines. We get del Pezzo surface
$X_1\subset \PP^6$ of degree $6$.

In the general algebraic case the surface $X_1$ may be not
isomorphic to $X$, i. e. $\Phi_{6,2}$ is not always a birational
involution (but I do not know such examples). However, in our
particular case $\Phi_{6,2}$ is a birational involution of Bertini
type. Indeed, if $\sigma\colon X^{\prime\prime}\rightarrow
X^\prime\rightarrow X$ is a blow up of $G$-orbit $\{Q_1,Q_2,Q_3\}$
together with $\{P_1,P_{-1}\}$, then $K_{X^{\prime\prime}}^2=1$.
But $-K_{X^{\prime\prime}}$ is not ample (because it is
degenerated del Pezzo surface of degree $1$ with three
$(-2)$-curves, or, after their contraction, three ordinary double
points). Nevertheless, the linear system
$|-2K_{X^{\prime\prime}}|$ gives contraction of $(-2)$-curves
composed with double covering $X^{\prime\prime}\rightarrow
Q^*\subset \PP^3$ of the quadratic cone branched over a fixed
point and a curve on $Q^*$ with three ordinary double points.

The involution of this double covering
$X^{\prime\prime}\rightarrow X$ induces the birational involution
$\Phi_{6,2}$. It differs from the standard Bertini involution
given by the Picard group, because $5$ blowing up points are not
in the general position. The point $P$ and ordinary double points
are still fixed under this involution. It also induces the
involution on the cubic surface $X_1$ that we consider in remark
$1.6$.

The birational involution $\Phi_{6,2}$ maps the pencil of rational
curves $\Pi=|-K_X-2P-P_1-P_{-1}|$ onto itself and changes the
components of its three fibers: $\Gamma_x\leftrightarrow
\Delta_x$, $\Gamma_y\leftrightarrow \Delta_y$,
$\Gamma_z\leftrightarrow \Delta_z$ (see $1.4.3$).

Apparently, we can write the particular equation for the residual
$G$-invariant pair of the curves of genus $0$ on $X$ (but we do
not need that). This $3$-linear $3$-homogenous equation should be
as follows.

\begin{multline*}
H\colon \sum_{g\in G} g((x_1-\lambda x_0)(y_1-\lambda
y_0)(z_1-\lambda^2z_0)+(x-\lambda x_0)(y_1-\lambda^2
y_0)(z_1-\lambda^2 z_0))=0.
\end{multline*}

The curve $H\sim -K_X$ is $G$-invariant, passes through $P_1$,
$P_{-1}$, and have the ordinary double singularities in them. This
means that it is reducible.

The link $\Phi_{6,2}\colon X\dashrightarrow X$ maps the linear
system $\mathcal H\sim -a K_X$ with the maximal singularity
$x=\{P_1, P_{-1}\}$ of multiplicity $r=r(x)$ to the linear system
$\mathcal H_1\sim -a_1 K_X$ with the basic cycle $x_1=x=\{P_1,
P_{-1}\}$ with the multiplicity $r_1 < a_1$ (this is not maximal
singularity for $\mathcal H_1$). The formulas for the coefficients
are given in~\cite{Isk3}, theorem 2.6, case $K_X^2=6$, case d):
\begin{equation}\label{8}
a_1=2a-r(x)<a,\ \ \ \ \ \ \ \ \ \ r_1=3a-2r(x)<a_1.
\end{equation}

So, $x=\{P_1,P_{-1}\}$ already is not the maximal singularity for
$\mathcal H_1$. The map $\chi_1=\chi\circ \Phi^{-1}\colon X
\dashrightarrow Y$ defined by $\mathcal H_1$ is not isomorphism,
so, by lemma--Noether inequality, $\mathcal H_1$ must have the
maximal singularity. Now it can be only the point $P$ of
multiplicity $r_1=r(P)>a_1$.

\subsubsection {\bf The case $d=1$.}
The only fixed point $x_1=P$. It does not lie on the
$(-1)$-curves, so it is in a general position. The corresponding
link $\Phi_{6,1}\colon X\dashrightarrow X_2$, where $X_2\subset
\PP^3$ is a smooth quadric with $\pic^G(X_2)=\frac{1}{2}\ZZ(-K_X)$
which we consider in section $1.7$, is the birational projection
$p_2\colon X\dashrightarrow X_2\subset \PP^3$ from the tangent
plane in $\PP^6$ to $X$ at $P$. This link blows up the point $P$
and contracts three conics $\Gamma_x$, $\Gamma_y$, $\Gamma_z$ to
the $G$-orbit $\mathcal A$ of length $3$ on the conic $C_0\subset
X_2$ which is the image of blowing up $P$. The stabilizator of $P\in
X$ is the whole $G$ which acts on the projectivisation of tangent
plane to $C_0$ as $S_3$ by the standard irreducible
$2$-dimensional linear representation. We completely considered
this situation in the lemma $1.8$.

There is $G$-invariant curve $C_0\sim -\frac{1}{2}K_{X_2}$, which
generates $\pic^G(X_2)$ on $X_2$. This means that, under the
notations from $1.8$, the formulas for the action of $\Phi_{6,1}$
are the following.
$$
-K_X\rightarrowtail -\frac{3}{2} K_{X_2}-2\mathcal A\sim
3C_0-2\mathcal A,
$$
\begin{equation}
\label{9} P\rightarrowtail -\frac{1}{2}K_{X_2}-\mathcal A\sim
C_0-\mathcal A,
\end{equation}
$$
\mathcal H_1\subset |-a_1K_{X}-r_1P|\rightarrowtail \mathcal
H_2\subset |-a_2K_{X_2}-r_2\mathcal A|, \ \ \ \ \mbox{\rm where }
$$
$$
a_2=\frac{3}{2}a_1-\frac{1}{2}r_1, \ \ \ \ r_2=2a_1-r_1.
$$
Under our hypothesis of maximality of the singularity in $P$ we
have $r_0=mult_P \mathcal H_1>a_1$, so $a_2<a_1$ and $r_2<a_2$.
This means that the $G$-orbit $\mathcal A$ is not the maximal
singularity for $\mathcal H_2$.

\subsection{} Now we are on the quadric $X_2\subset \PP^3$ with
linear system $\mathcal H_2\sim -a_2 K_{X_2}$. Denote this system
by $\mathcal H$ and put $a=a_2$ for simplicity. The map
$\chi_2\colon X_2\dashrightarrow Y=\PP^2$ is not isomorphism, so,
by Noether inequality, there is a maximal singularity, $G$-orbit
$x$ of multiplicity $r=r(x)>a$ and of length $d<K_{X_2}^2=8$.
$d|12=|G|$so the only possibilities for $d$ are $d=1,2,3,4,6$. We
can exclude the cases $d=1$ and $d=4$ by lemma $1.8$. So consider
the other cases $d=2,3,6$.

\subsubsection{} {\bf The case $d=2$}.
There is two $G$-orbits $\{R_1,R_2\}$ and $\{P_1,P_{-1}\}$ (see
$1.8$). Start from $x=\{R_1,R_2\}=C_0\cap C_1$ of multiplicity
$r=r(x)>a$. Such $R_1$, $R_2$ do not lie on one line on the
quadric $X_2$, because each line intersects $C_0$ by a unique
point ($C_0$ is the hyperplane section). So, $R_1$ and $R_2$ are
in the general position. There exists link $\Phi_{8,2}$
(see~\cite{Isk3}, $2.6$, the case $K_X^2=8$, a), $d=2$), which is
a blow up $X_2\dashrightarrow X_1$ of the points $\{R_1,R_2\}$.
There is $G$-invariant structure of conic bundle $\pi\colon
X_1\rightarrow \PP^1$ on $X_1$, which is given by the pencil
$\Pi_0$ from $1.8$ (iv). Let $\mathcal H_1$ be the strict
transform on $X_1$ of linear system $\mathcal H$. Then $\mathcal
H_1\sim -a_1K_{X_1}+b_1f_1$, where $f_1\cong C_0^\prime$ is
$G$-invariant fiber. We have
$$
K_{X_1}^2=6, \ \ \ \ -K_{X_1}=E+C_0^\prime+C_1^\prime \sim E+2f_1,
$$
where $E$ is the exceptional divisor on $X_1\rightarrow X_2$, i.
e. $G$-invariant pair of $(-1)$-curves and $C_0^\prime$,
$C_1^\prime$ are strict transforms of $C_0$, $C_1$.
By~\cite{Isk3}, $2.6$
\begin{equation}
\label{10}
a_1=2a-r(x)<a,
\end{equation}
$$
b_1=2(r(x)-a).
$$

From~\ref{10} one can see that $a_1<a$ and $b_1>0$. This means
that by the Noether inequality (because $X_1\ncong Y$) there is a
maximal singularity in $\mathcal H_1$. But:

\subsubsection \Lem {\it There is no zero-dimensional $G$-orbits $x_1$
of length $d_1$ on $X_1$, whose points are in the general position
in the sense of $2.1$ {\rm 2)}, i. e. there is no point in it that
lies on the degenerated fiber and there are no two points that lie
on the same fiber of morphism $\pi\colon X_1\rightarrow \PP^1$
{\rm (}this is a condition of existence of  ``untwisting'' link
for $\mathcal H_1$ on the conic bundle with $b_1\geq 0${\rm )}. }
\medskip

\Proof. By lemma $1.8$ (iv), a cyclic group of order $2$ acts on
$\PP^1$, so for generality the length $d_1$ must be less or equal than
$2$. There is no $G$-fixed points on $X_1$ (as on $X_2$) The orbit
$\{P_1,P_{-1}\}$ of length $2$ by lemma $1.8$ lies on the
degenerated fibers. On the non-degenerated fibers acts $S_3$
without fixed points. \qed

\medskip

So, the decomposition by the algorithm in this situation is
impossible. This means that $\{R_1,R_2\}$ is not maximal, and we
need to come back to the quadric $X_2$. Consider the other
possibilities.

\subsubsection{} {\bf The case $d=2$} with maximal singularity
$x=\{P_1,P_{-1}\}$, $r(x)>a$. As in $2.4.1$, the pair
$\{P_1,P_{-1}\}$ is in the general position and there exists link
$\Phi_{8,2}$ (see~\cite{Isk3}, theorem 2.6, the case $K_X^2=8$, c),
$d=2$), which is the blow up $X\dashrightarrow X_1$ of $P_1$,
$P_{-1}$. There is a structure of conic bundle $\pi\colon
X_1\rightarrow \PP^1$ on $X_1$ that is given by the pencil $\Pi_1$
from $1.8$ (v).

The difference from the previous case is the following. The action
of $G=S_3\times \ZZ_2$ on $\pi\colon X_1\rightarrow \PP^1$ is
given by two commutating actions. $S_3$ acts on the base $\PP^1$
without fixed points and this action can be shifted to $X_1$ as
the changing of the fibers. $\ZZ_2$ acts on $X_1\rightarrow \PP^1$
only on the fibers. The fixed points in this case are the
intersection of the fibers with the curve $C^\prime_0\subset X_1$,
which is the pre-image of the curve $C_0\subset X_2$. So, $\ZZ_2$
acts as the classical De Jonquieres involution. The fixed curve
$C_0^\prime$ determines $\pi\colon X_1\rightarrow \PP^1$. Indeed,
the factor by this involution is the ruled surface
$\FF_1\rightarrow \PP^1$ branched over $C_0^\prime$. The pair of
the $(-1)$-curves $E$ covers the exceptional section $S_1\subset
\FF_1$.

Notice that there is no fixed points on the base $\PP^1$ (so there
is no $G$-invariant fiber as in the previous case). But the class
of the fiber $[f_1]$ in $\pic^G(X)$ is $G$-invariant, so the
action of this link is the same as in (\ref{10}).

By the Noether inequality there must be the maximal singularity in
the linear system $\mathcal H\sim -a_1K_{X_1}+b_1f_1$ which is
$G$-orbit $x_1\subset X_1$ of multiplicity $r_1=r(x_1)>a_1$, whose
points are in the general position in the sense of $2.4.2$.

There is no $G$-fixed points on $X_1$ and the only $G$-invariant
orbit of length $2$, $\{R_1,R_2\}$, lies in the degenerated fibers
(see $1.8$ (v)).

The both $G$-orbits of length $3$, $\mathcal A$ and $\mathcal B$,
satisfy the generality condition. Both of them lie on $C_0^\prime$
and they change under the involution of double covering
$\pi|_{C_0^\prime}\colon C^\prime_0\rightarrow \PP^1$ (do not
confuse with involution of the action of $\ZZ_2$).

There are also $G$-orbits of length $6$. By the generality
condition they lie also on $C_0^\prime$, one point on the fiber.
The other intersection points of $C_0^\prime$ with fibers are
complementary $G$-orbit of length $6$.

There is an ``untwisting'' link concerned with every maximal
singularity $x_1\subset X_1$ of length $3$ or $6$. That is, the
elementary transform
\begin{equation}
\label{11} \xymatrix{
&\Phi\colon X_1\ar@{->}[rd]\ar@{-->}[rr]^{\mathcal E_{x_1}}&&X_1^\prime.\ar@{->}[ld]& \\
&&\PP^1 }
\end{equation}

Remind that the elementary transform blows up the points of orbit
$x_1\subset X_1$ and contracts the strict transforms of the fibers
on witch these points lie. In the general case an elementary
transform is not a birational involution. But in our case it is.
Indeed, the transform $\mathcal E_{X_1}$ induces the isomorphism
on the curve $C_0^\prime$, and, as we notice above, $C_0^\prime$
completely determines $\pi\colon X_1\rightarrow \PP^1$, so
$X_1^\prime=X_1$ in (\ref{11}).

The birational involution $\Phi=\mathcal E_x$ acts by the
following formulas from~\cite{Isk3}, $2.6$.
$$
-K_{X_1}\rightarrowtail -K_{X_1}+d_1f_1-2x_1,\ \ \ \ a_1^\prime
=a_1,
$$
$$
\ \ \ \ \ \ \ \ \ \ \ \ \ \ \ \ \ f_1\rightarrowtail f_1,\ \ \ \ \
\ \ \ \ \ \ \ \ \ \ \ \ \ \ \ \ \ \ \ \ \ \ \ \ \ \ \ b_1^\prime =
b_1+d_1(a_1-r_1),
$$
$$
\ \ \ \ \ \ \ \ \ \ x_1\rightarrowtail d_1f_1-x_1,\ \ \ \ \ \ \ \ \
\ \ \ \ \ \ \ \ \ \ \ \ \ \ \ r_1^\prime = 2(r_1-a_1),
$$
where $d_1$ is the length of the orbit $x_1$. As the last formula
shows, the multiplicity of the maximal singularities decreases.
After the finite number of such transforms we can ``untwist'' all
maximal singularities (which means that the multiplicities of the
base points in the modified linear system will be not greater then
$a$). Then, by the Noether inequality, the coefficient at $b_1$ at
the fiber will be less than $0$.

In this situation the link $\Phi^{-1}_{8,2}\colon X_1\rightarrow
X_2$ (the inverse to $\Phi_{8,2}$) contracts $G$-invariant pair of
$(-1)$-curves $E$ (the sections of the pencil $\pi\colon
X_1\rightarrow \PP^1$) to $G$-orbit of length two
$x=\{P_1,P_{-1}\}$. It acts by the following formulas.
\begin{equation}
\label{12} -K_{X_1}\rightarrowtail -K_{X_2}-x, \ \ \ \
a=a_1+\frac{2}{3}b_1,
\end{equation}
$$
\ \ 2f_1\rightarrowtail -K_{X_2}-2x, \ \ \ \ \ \ r(x)=a_1+b_1.
$$
These formulas involves $a<a_1$ and $r(x)<a$, because $b_1<0$.
This means that this link is untwisted and $x=\{P_1, P_{-1}\}$
already is not the maximal singularity.

\subsection{}
The last cases are those whose singularities of $\mathcal H$ are
$G$-orbits of length $3$ or $6$. There are only two orbits of
length $3$: $\mathcal A$ and $\mathcal B$ (see $1.8$). First
consider

\subsubsection{} {\bf The case $d=6$}. Suppose that there is a maximal
singularity, $G$-orbit $x$ of length $6$, all whose points are in
the general position and the multiplicity $r=r(x)>a$, where
$\mathcal H\sim -aK_{X_2}$.

Then there is an untwisting link $\Phi_{8,6}$, the Geiser's
involution (see~\cite{Isk3}, $2.6$). Indeed, if we blow up $x$,
then we get del Pezzo surface of degree $2$, on which the
classical Geiser's involution acts biregular. Apply this
involution and contract the blowing up divisor. The link
$\Phi_{8,6}$ acts by the following formulas.
$$
a_1=7a-6r(x)<a,
$$
$$
r(x_1)=8a-7r(x)<a_1.
$$

\subsubsection{} {\bf The case $d=3$}. If the maximal singularity
is $x=\mathcal A$, then the link $\Phi_{8,3}$ with center in this
singularity is the inverse map $\Phi_{6,3}^{-1}\colon
X_2\dashrightarrow X$ and we are in the situation from which we
start. If the maximal singularity is $x=\mathcal B$ (we missed
this case in~\cite{Isk1}), then the link $\Phi_{8,3}^\prime\colon
X_2 \dashrightarrow X^\prime$ acts from $X_2$ to del Pezzo surface
of degree $6$ $X^\prime$. But $X^\prime \cong X$, because the
birational involution on $X$ with the center in $\{P_1,P_{-1}\}$
that we consider in $2.3.2$ changes $\{\Gamma_2, \Gamma_y,\Gamma_z
\} \leftrightarrow \{\Delta_x,\Delta_y,\Delta_z\}$. Thus, it
changes the cycles $\mathcal A$ and $\mathcal B$ on $C_0$ that
correspond to these directions in the tangent plane to $P\in X$.

So, we consider all possibilities but have not found any
$G$-equivariant $G$-map $X\dashrightarrow Y=\PP^2$. So there is no
such map, which prove the
proposition~\ref{Proposition:main_result}.

\subsection{Commentary}
In the proof we determine that the smooth relative minimal
$G$-models on which there exists $G$-equivariant map from $X$ are
only the quadric $X_2$ or the conic bundle with two different
actions of $G$. It is interesting to find by this method all of
$G$-surfaces on which $Y=\PP^2$ may be $G$-equivariant mapped.

It is also interesting the following. If we consider only $S_3$
instead of $G$ with the similar actions on $X$ and $Y$, then these
surfaces are $S_3$-birational equivalent. Indeed, in this case
there are three fixed points $P$, $P_1$ and $P_{-1}$ on $X$.
Therefore, there is two fixed points $P_1$, $P_{-1}$ on the
quadric $X_2$. The stereographic projection from one of them is
the required $G$-equivariant birational isomorphism, so both of
these embeddings of $S_3$ it $Cr_2(\CC)$ are conjugate. From the
other hand, the involution $(x,y)\rightarrow (x^{-1},y^{-1})$ of
course is conjugate to the linear one.

\end{document}